\documentclass[10pt,a4paper]{article}
\usepackage{amsfonts}
\usepackage{amssymb}
\usepackage{mathrsfs}
\usepackage{amsthm}
\usepackage{setspace}
\usepackage{color}
\usepackage{stmaryrd}

\usepackage{fancyhdr}
\newtheorem{theo}{Theorem}[section]
\newtheorem{prop}[theo]{Proposition}
\newtheorem{cor}[theo]{Corollary}
\newtheorem{lemm}[theo]{Lemma}

\theoremstyle{definition}

\newtheorem{defi}[theo]{Definition}

\newenvironment{lproof}{\emph{Proof of Lemma.}}{ \qed \par}

\newcommand{\be}{\begin{eqnarray*}}
\newcommand{\ee}{\end{eqnarray*}}
\newcommand{\beqa}{\begin{eqnarray}}
\newcommand{\eeqa}{\end{eqnarray}}
\newcommand{\ba}{\begin{array}}
\newcommand{\ea}{\end{array}}

\newcommand{\mf}{\mathfrak}

\newcommand{\mbb}{\mathbb}

\begin{document}

\title{Reducing almost Lagrangian structures and almost CR geometries to partially integrable structures}
\author{Stuart Armstrong}
\date{July 2008}
\maketitle

\begin{abstract}
This paper demostrates a method for analysing almost CR geometries $(H,J)$, by uniquely defining a partially integrable structure $(H,K)$ from the same data. Thus two almost CR geometries $(H,J)$ and $(H',J')$ are equivalent if and and only if they generate isomorphic induced partially integrable CR geometries $(H,K)$ and $(H',K')$, and the set of CR morphisms between these spaces contains an element that maps $J$ to $J'$. Similar results hold for almost Lagrangian structures.
\end{abstract}

\section{Introduction}
CR geometry is a particularly rich mathematical seam, spawning elegant results and successful applications all over the place with joyful abandon. Though partially integrable CR manifolds of hyper-surface type can be sucessfully analysed by constructing a unique normal Cartan connection (see \cite{CR} and \cite{CartEquiv}), the same is not true for general almost CR structures of the same type, where no such constructions exist.

However, this paper details a way of locally matching each `generic' almost CR structure $(H,J)$, in a purely algebra\"ic fashion, to a uniquely defined partially integrable CR structure $(H,K)$. The data of $(H,J)$ can thus be fully captured by the normal Cartan connection for $(H,K)$ and by $J$. In particular, two almost CR structures are equivalent if and only if they generate isomorphic induced partially integrable CR structures $(H,K)$ and $(H',K')$, and the set of CR morphisms between these spaces contains an element that maps $J$ to $J'$.

Similar results hold for Lagrangian geometries, which are another real form of these CR geometries.

The construction of $K$ from $J$ requires certain assumptions on the eigenvalues of an automorphism $A$ of $H$, specifically that they all be non-real. This is the `generic' condition needed for the almost CR structure -- note that the condition is vacuous in the definite signature case, where $A$ has only imaginary eigenvalues. When the condition does fail, the structure on the manifold can best be described as a mixed structure, intertwining Lagrangian and CR structures.

\subsection*{Acknowledgements}
It gives me great pleasure to acknowledge the financial support of project P19500-N13 of the ``Fonds zur F\"orderung der wissenschaftlichen Forschung (FWF)''.

\section{The CR case}

An almost CR structure of hypersurface type is given by:
\begin{itemize}
\item a manifold $M$ of dimension $1+2n$,
\item a distribution $H \subset TM$ of rank $2n$ generating a ``contact'' form $\omega \in \Gamma(\wedge^2 H^* \otimes (TM/H))$, which may be degenerate,
\item an almost-complex structure $J$ on $H$, such that $\nu$ is non-degenerate, where $\nu$ is defined as
\be
\nu(X,Y) = \frac{1}{2} \big((\omega(X,Y) + \omega(JX,JY) \big),
\ee
for any sections $X$ and $Y$ of $H$.
\end{itemize}
Note that $\nu$ is $J$-hermitian, in that $\nu(JX,JY) = \nu(X,Y)$ for all $X,Y \in \Gamma(H)$. It defines a metric $g$ given by
\be
g(X,Y) = \nu(X,JY).
\ee
The signature of an almost CR structure is the signature of $g$. If $g$ is positive definite, then the non-degeneracy of $\nu$ implies the non-degeneracy of $\omega$; this is not the case for other signatures.
\begin{defi}
A partially integrable CR manifold is one defined as above with $\nu = \omega$.
\end{defi}

We may construct an automorphism $A$ of $H$ as $A = g^{-1} \omega$. Then the core theorem of this paper is:
\begin{theo} \label{main:theo}
If $\omega$ and $\nu$ are both non-degenerate, and none of the eigenvalues of $A$ is real, then there exists an almost-complex structure $K$ on $H$, uniquely defined by the data $(H,J)$, such that $(H,K)$ defines a partially integrable CR manifold.
\end{theo}
In the definite signature case, all the eigenvalues of $A$ must be pure imaginary, so the restriction on $A$ is not needed. It is an open condition that will be satisfied by ``most'' distributions, and certainly those obtained by small deformations of partially integrable CR structures: for in that case $A = J$, with eigenvalues $\pm i$.

Since partially integrable CR structures are defined by a unique normal Cartan connection (\cite{CR} and \cite{CartEquiv}), this has the immediate corollary that:
\begin{cor}
Two almost CR structures $(H,J)$ and $(H',J')$ are isomorphic if and only if their generated partially integrable CR structures $(H,K)$ and $(H',K')$ define isomorphic normal Cartan connections, and if there exists a CR morphism between them that pulls back $J$ to $J'$.
\end{cor}

The methods used to prove Theorem \ref{main:theo} involve constructing $K$ linearly from $J$ at each point of $M$. The procedure is easily seen to be continuous, generating a continuous $K$. Now pick any point $x$ in $M$, and let $V = H_x$; by an abuse of notation, we will drop the indexes and refer to $\omega_x$, $\nu_x$ $g_x$, $A_x$ and $J_x$ as $\omega$, $\nu$, $g$, $A$ and $J$.

The endomorphism space $V\otimes V^*$ decomposes under the triple $(\nu,g,J)$ as
\be
V \otimes V^* = \ \mf{su}(g,J) \ \oplus \ (\wedge^{1,1}_0 V) \ \oplus \ (\wedge^2_{\mbb{C}} V) \ \oplus \ (\odot^2_{\mbb{C}} V) \ \oplus \ (\mbb{R} \cdot J) \ \oplus \ (\mbb{R} \cdot Id).
\ee

Relevant subalgebras of this are the complex algebra $\mf{sl}(V,J) = \mf{su}(g,J) \oplus (\wedge^{1,1}_0 V)$, the symplectic algebra $\mf{sp}(\nu) = \mf{su}(g,J) \oplus (\odot^2_{\mbb{C}} V) \oplus (\mbb{R} \cdot J)$, the orthogonal algebra $\mf{so}(g) = \mf{su}(g,J) \oplus (\wedge^2_{\mbb{C}} V) \oplus (\mbb{R} \cdot J)$ and the conformal unitary algebra $\mf{su}(g,J) \oplus  (\mbb{R} \cdot J) \oplus  (\mbb{R} \cdot Id)$. Following the tradition of parabolic geometry (\cite{TBIPG}, \cite{TCPG}), this last algebra will be designated $\mf{g}_0$ and the corresponding group $G_0$. It is the largest group that preserves a given CR structure. Upon constructing $K$ from $J$, we will see that the definition is unique up to $G_0$ action, hence that it defines a unique CR structure.

Note that $\mf{so}(g) \cap \mf{gl}(V,J) = \mf{g}_0$, hence that $SO(g) \cap GL(V,J) = G_0$. The process for construction $K$ flows from:
\begin{prop}
For any given collection of non-degenerate $\omega$, $J$, $\nu$ $A$ and $g$, defined as above with $A$ having only non-real eigenvalues, there exists an element $e$ in $SO(g)$ such that $\omega$ is hermitian for the almost-complex structure $K = e^{-1} J e$. This $e$ is defined up to the left action of $G_0$; since $G_0$ preserves $J$, this defines $K$ uniquely.
\end{prop}

Constructing this $e$ is basic linear algebra. We will be operating in the complexified space $V_{\mbb{C}} = V \otimes \mbb{C}$, keeping track of the subspace $V = V \otimes 1$ via complex conjugation.

\begin{defi}
The space $D_{\alpha}^k$ is defined to be the $k$-th generalised eigenspace for $A$ with eigenvalue $\alpha$, i.e.~the kernel of the linear map $(A - \alpha Id)^k$. The Jordan normal form decomposition of $A$ demonstrates that
\be
V_{\mbb{C}} = \oplus_j D_{\alpha_j},
\ee
where $\alpha_j$ are the eigenvalues of $A$ and $D_{\alpha_j} = D_{\alpha_k}^k$ for some $k$ where $D^{k+1}_{\alpha_j} = D^{k}_{\alpha_j}$.
\end{defi}

\begin{defi}
The set $S$ is defined to be the set of eigenvalues of $A$; since $A$ is non-degenerate, $0 \notin S$, and by assumption, $S \cap \mbb{R} = \emptyset$.
\end{defi}

\begin{lemm}
The space $D_{\alpha}$ is $g$-orthogonal to all eigenspaces $D_{\beta}$, except when $\beta = -\alpha$. Moreover, $g$ gives a non-degenerate pairing between $D_{\alpha}$ and $D_{-\alpha}$.
\end{lemm}
\begin{lproof}
Let $u \in D_{\alpha}^1$, $v \in D_{\beta}^1$. Then since $g(Au,v) = -g(Av,u)$, we must have $\alpha g(u,v) = -\beta g(u,v)$. Hence either $g(u,v) = 0$, or $\alpha = -\beta$. Assume for the moment that $\alpha \neq - \beta$; hence $D_{\alpha}^1 \perp D_{\beta}^1$.

Reasoning by induction, assume that $D_{\alpha}^j \perp_g D_{\beta}^k$, and let $u \in D_{\alpha}^{j}$, $v \in D_{\beta}^{k+1}$. Then
\be
\alpha g(u,v) = g(Au,v) = = -g(u,Av) = -g(u,\beta v + v^k) = -g(u,\beta v).
\ee
Hence $D_{\alpha}^j \perp_g D_{\beta}^{k+1}$. Since we may induct both $j$ and $k$, and since the generalised eigenspaces stabilise after finitely many steps, we must have $D_{\alpha} \perp D_{\beta}$.

Consequently, the pairing under $g$ gives $D_{\alpha} \supset D_{-\alpha}^*$ and $D_{-\alpha} \supset D_{\alpha}^*$. Dimensional considerations imply that $g$ pairs $D_{\alpha}$ and $D_{-\alpha}$ in a non-degenerate fashion.
\end{lproof}
Since $A$ is real, it commutes with complex conjugation, implying that $\overline{S} = S$ and that $\overline{D_{\alpha}} = D_{\overline{\alpha}}$.

Define $S_+$ as the set of elements $\alpha$ in $S$ such that the argument of $\alpha$ is in $(0, \pi/2]$. Then $S = S_+ \cup \overline{S_+} \cup -S_+ \cup -\overline{S_+}$. For any $\alpha$ in $S_+$, define the space
\be
C_{\alpha} = D_{\alpha} + D_{\overline{\alpha}} + D_{-\alpha} + D_{-\overline{\alpha}}.
\ee
These $C_{\alpha}$ are mutually orthogonal, non-degenerate under $g$ and closed under complex conjugation. This means that $C_{\alpha} = C_{\alpha}^{\mbb{R}} \otimes \mbb{C}$, where $C_{\alpha}^{\mbb{R}}$ is a real subspace of $V$. If $\alpha$ is not pure imaginary, then $C_{\alpha}^{\mbb{R}}$ must be of split signature $(2p,2p)$, since $D_{\alpha} + D_{\overline{\alpha}}$ must be the complexification of an even dimensional isotropic space, pairing non-degenerately with $D_{-\alpha} + D_{-\overline{\alpha}}$.

If $\alpha$ is pure imaginary, let $v + \overline{v}$ be an orthonormal element of $C_{\alpha}^{\mbb{R}}$, for $v \in D_{\alpha}$. Then
\be
||iv -i\overline{v}||^2 &=& -2 g(iv,i\overline{v}) + g(iv,iv) + g(\overline{v},\overline{v})\\
&=& 2g(v,\overline{v}) \\
&=& ||v + \overline{v}||^2,
\ee
since $v$ and $\overline{v}$ are isotropic. Hence $C_{\alpha}^{\mbb{R}}$ is of signature $(2p,2q)$.

Let $L_{\pm}$ be the $\pm i$ eigenspace of $J$. These two spaces must be isotropic with $L_+ \oplus L_- = V$, by the properties of $J$. Note that $\overline{L_+} = L_-$. By the signature results for $C_{\alpha}^{\mbb{R}}$, we have the following lemma:
\begin{lemm}
There exists a decomposition of $L_+$ as
\be
L_+ = \oplus_{\alpha \in S_+} P_{\alpha},
\ee
such that the spaces $Q_{\alpha} = P_{\alpha} \oplus \overline{P_{\alpha}}$ are mutually orthogonal, non-degenerate under $g$, closed under complex conjugation, and of same dimension and signatures as $C_{\alpha}$.
\end{lemm}
\begin{lproof}
Pick an orthonormal basis in $V$ for the hermitian metric $g + i \nu$, group the basis elements together to generate subspaces of the correct signature, and complexify into subspaces of $V_{\mbb{C}}$. Since these spaces are all preserved by $J$, they give the required splitting of $L_+$.
\end{lproof}

We now choose a map $e$ on $V_{\mbb{C}}$, defined in the following way: for $\alpha$ not purely imaginary, let $P_{\alpha} = P^1_{\alpha} \oplus P^2_{\alpha}$, where $P^j_{\alpha} \oplus \overline{P^j_{\alpha}}$ is isotropic. Then map $D_{\alpha}$ into $P^1_{\alpha}$ in any fashion, map $D_{\overline{\alpha}}$ into $\overline{P^1_{\alpha}}$ by conjugation, $D_{-\alpha}$ into $\overline{P^2_{\alpha}}$ by duality under $g$, and $D_{-\overline{\alpha}}$ into $P^2_{\alpha}$ by duality then conjugation (or conjugation then duality -- $g$ is real, hence commutes with conjugation).

For $\alpha$ purely imaginary, pick an orthonormal basis $v_j \oplus \overline{v_j}$ of $C^{\mbb{R}}_{\alpha}$ with $v_j \in D_{\alpha}$, an orthonormal basis $u_j \oplus \overline{u_j}$ of $Q^{\mbb{R}}_{\alpha}$ for $u \in L_+$, and map one basis to the other, sending $D_{\alpha}$ into $P_{\alpha}$.

By construction, $e$ preserves the metric $g$ and complex conjugation; thus it is an element of $SO(g)$. Let $e'$ be another element of $SO(g)$ that maps $D_{\alpha}$ to $L_+$ whenever $\alpha$ has positive imaginary part; then $e' = fe$, where $f$ is an element of:
\be
(GL(L_+) \oplus GL(L_-)) \cap SO(g).
\ee
But this intersection is $G_0$, as the real part of $(GL(L_+) \oplus GL(L_-))$ is just $GL(V,J)$. Thus $e$ is unique up to left $G_0$ action. By construction, the endomorphism $eAe^{-1}$ must have eigenspaces that are subspaces of $L_{\pm}$, and hence commute with $J$. Since $e$ is orthogonal, this implies that $e^{-1} \cdot \omega$ is $J$-hermitian, where
\be
(e^{-1} \cdot \omega)(X,Y) = \omega(e^{-1}X,e^{-1}Y).
\ee
\begin{lemm}
The form $\omega$ is hermitian under the complex structure $K = e^{-1} J e$, i.e.
\be
\omega(KX,KY) = \omega(X,Y),
\ee
and $K$ is invariantly defined independently of the choice of $e$.
\end{lemm}
\begin{lproof}
First note that
\be
\omega(KX,KY) &=& \omega(e^{-1}J(eX),e^{-1}J(eY)) \\
&=& (e^{-1} \cdot \omega) (J(eX),J(eY)) \\
&=& (e^{-1} \cdot \omega) (eX,eY) \\
&=& \omega(e^{-1}(eX),e^{-1}(eY)) \\
&=& \omega(X,Y),
\ee
since $(e^{-1} \cdot \omega)$ is $J$-hermitian. Now let $e' = fe$ be another suitable map. Then
\be
(e')^{-1} J e' = e^{-1} f^{-1}Jfe = e^{-1} Je = K,
\ee
since $f \in G_0$ preserves $J$.
\end{lproof}

\section{The Lagrangian case}

Almost Lagrangian structures are given by a distribution $H$ with contact form $\omega$, as above, and by a decomposition
\be
H = E \oplus F,
\ee
into two bundles of equal size. This can be characterise by the existence of a trace-free involution $\sigma$ squaring to the identity, with $E$ as its $+1$ eigenspace and $F$ its $-1$ eigenspace. Partial integrability is given by the relation:
\be
\omega(\sigma(X),\sigma(Y)) = -\omega(X,Y),
\ee
equivalent to the isotropy of $E$ and $F$ under $\omega$ (notice the change in sign compared with the CR case). The canonical two-form that we will need is $\nu$, defined by
\be
\nu(X,Y) = \frac{1}{2} \big( \omega(X,Y) - \omega(\sigma(X),\sigma(Y)) \big),
\ee
while the (split) metric $g$ is
\be
g(X,Y) = \nu(X,\sigma(Y)).
\ee
If $\omega$ and $\nu$ are non-degenerate, and the automorphism $A = g^{-1} \omega$ does not have any purely \emph{imaginary} eigenvalues, then the proof proceed as in the CR case (except that now $D_{\alpha} \oplus D_{\overline{\alpha}}$ will be mapped into $L_+$, rather than $D_{\alpha} \oplus D_{-\overline{\alpha}}$ as was the case then; note also that $\overline{L}_{\pm} = L_{\pm}$, $L_+ = E \otimes \mbb{C}$ and $L_- = F \otimes \mbb{C}$).

\section{Real and imaginary eigenvalues}
If $A$ in the CR case has a real eigenvalue, the above procedure does not work. Since a $g$-skew automorphism with real eigenvalues can be approximated arbitrarily closely by those with complex eigenvalues, it must still remain the case that $D_{\alpha}$ is even-dimensional, and that $C_{\alpha} = D_{\alpha} \oplus D_{-\alpha}$ is of split signature $(2p,2p)$. The only obstruction to choosing $e$ as above is that $e$ cannot now be chosen to lie inside $SO(g)$, but only in its complexification $SO(g,\mbb{C})$ (uniqueness of the definition of $e$ is preserved by assigning $D_{\alpha}$ to $L_+$ when $\alpha >0$). However, a different approach can bear fruit.

Define $M_A$ to be the subset of $M$ where $A$ has real eigenvalues. It must be closed, implying that $M^c = M - M_A$ is open, hence that $M^c$ is a submanifold of $M$. On $M^c$, we have a unique choice of $K$, but this choice cannot necessarily be extended continuously across $M_A$. If $M_A$ has empty interior, and does not separate $M$ into components, then the degeneracy on $M_A$ does not matter much: the equivalence problem for $(H,J)$ can be analysed away from $M_A$, and extended to $M_A$ by continuity. Even if $M$ does get separated in to components, the equivalence problem can still be analysed on each component separately, and the resulting limits ``glued together''.

If $M_A$ does have a non-empty interior, then the natural structure on it is an intertwined structure: a decomposition of $H$ into isotropic $H_1 \oplus H_2$, such that there exists a $J$ on $_1$ and a $\sigma$ on $H_2$ with $\omega|_{\wedge^2 H_1}$ being $J$-hermitian, and $\omega|_{\wedge^2 H_2}$ being $\sigma$-Lagrangian. On any connected subset of $M_A$ where the rank of the generalised eigenspaces for real eigenvalues is constant, such a structure can be defined. Since that rank is upper semi continous, and bounded above by $2n$, this allows us to partition $M_A$ into components where such structures are defined, excluding only sets of small dimension.

Of course, the converse results hold for almost Lagrangian structures with an $A$ with pure imaginary eigenvalues.
\bibliographystyle{alpha}
\bibliography{ref}

\end{document}